\documentclass[12pt,fullpage]{article}

\usepackage{amsfonts,graphicx,amsmath,xr,epsfig,comment}
\input epsf

\def\B{\mathbb{B}}

\def\R{\mathbb{R}}
\def\N{\mathbb{N}}
\def\epsilon{\varepsilon}

\newcommand{\be}{\begin{equation}}
\newcommand{\ee}{\end{equation}}
\newcommand{\baa}{\begin{array}}
\newcommand{\eaa}{\end{array}}
\newcommand{\ba}{\begin{eqnarray}}
\newcommand{\ea}{\end{eqnarray}}

\newtheorem{theo}{\bf Theorem}[section]
\newtheorem{lem}[theo]{\bf Lemma}

\newtheorem{rem}[theo]{\bf Remark}

\begin{document}
\date{}
\title{\bf{Asymptotic behavior of solutions of semilinear elliptic equations in unbounded domains: two approaches}}
\author{Messoud Efendiev$^{\hbox{\small{ a}}}$ and Fran\c cois Hamel$^{\hbox{\small{ a,b}}}$\thanks{The second author is indebted to the Alexander von~Humboldt Foundation for its support. He is also supported by the ANR project PREFERED.}\\
\\
\footnotesize{$^{\hbox{a }}$Helmholtz Zentrum M\"unchen, Institut f\"ur Biomathematik und Biometrie}\\
\footnotesize{Ingolst\"adter Landstrasse 1, D-85764 Neuherberg, Germany}\\
\footnotesize{$^{\hbox{b }}$Aix-Marseille Universit\'e \& Institut Universitaire de France}\\
\footnotesize{LATP, Avenue Escadrille Normandie-Niemen, F-13397 Marseille Cedex 20, France}}
\maketitle

\begin{abstract}
In this paper, we study the asymptotic behavior as $x_1\to+\infty$ of solutions of semilinear elliptic equations in quarter- or half-spaces, for which the value at $x_1=0$ is given. We prove the uniqueness and characterize the one-dimensional or constant profile of the solutions at infinity. To do so, we use two different approaches. The first one is a pure PDE approach and it is based on the maximum principle, the sliding method and some new Liouville type results for elliptic equations in the half-space or in the whole space~$\R^N$. The second one is based on the theory of dynamical systems.
\end{abstract}


\section{Introduction}\label{intro}

Let $\Omega$ be the domain of $\R^N$ ($N\ge 2$) defined by
$$\baa{rcl}
\Omega & = & (0,+\infty)\times\R^{N-2}\times(0,+\infty)\vspace{3pt}\\
& = & \big\{x=(x_1,x',x_N)\in\R^N\ \big|\ x_1>0,\ x'=(x_2,\ldots,x_{N-1})\in\R^{N-2},\ x_N>0\big\}.\eaa$$
This paper is devoted to the study of the large space behavior, that is as $x_1\to+\infty$, of the nonnegative bounded classical solutions $u$ of the equation
\be\label{eq}\left\{\baa{rcll}
\Delta u+f(u) & = & 0 & \hbox{ in }\Omega,\vspace{3pt}\\
u(x_1,x',0) & = & 0 & \hbox{ for all }x_1>0\hbox{ and }x'\in\R^{N-2},\vspace{3pt}\\
u(0,x',x_N) & = & u_0(x',x_N) & \hbox{ for all }x'\in\R^{N-2}\hbox{ and }x_N>0,\eaa\right.
\ee
where the function
$$u_0:\R^{N-2}\times(0,+\infty)\to\R_+=[0,+\infty)$$
is given, continuous and bounded. The solutions $u$ are understood to be bounded, of class~$C^2(\Omega)$ and to be continuous on $\overline{\Omega}\ \backslash\ \big(\{0\}\times\R^{N-2}\times\{0\}\big)$. From standard elliptic estimates, they are then automatically of class $C^{2,\beta}_b([\epsilon,+\infty)\times\R^{N-2}\times\R_+)$ for all $\epsilon>0$ and $\beta\in[0,1)$. Here and below, for any closed set $F\subset\R^N$ and $\beta\in[0,1)$, we write
\be\label{defC2betab}
C^{2,\beta}_b(F):=\Big\{u:F\to\R\ \big|\ \|u\|_{C_b^{2+\beta}(F)}=\sup_{x\in F} \|u\|_{C^{2,\beta}\big(\overline{B(x,1)}\cap F\big)}<\infty\Big\},
\ee
where $B(x,1)$ means the open euclidean ball of radius $1$ centered at $x$. Problems sets in the half-space
$$\Omega'=(0,+\infty)\times\R^{N-1}$$
will also be considered in this paper, see (\ref{eqbis}) below.

The value of $u$ at $x_1=0$ is then given and the goal is to describe the limiting profiles of~$u$ as $x_1\to+\infty$. If the equation were parabolic in the variable~$x_1$, we would then be reduced to characterize the $\omega$-limit set of the initial condition $u_0$. However, problem~(\ref{eq}) is an elliptic equation in all variables, including $x_1$, and the ``Cauchy" problem~(\ref{eq}) with the ``initial value"~$u_0$ at $x_1=0$ is ill-posed. There might indeed be several solutions~$u$ with the same value $u_0$ at $x_1=0$. Nevertheless, under some assumptions on the nonlinearity~$f$, we will see that the behavior as $x_1\to+\infty$ of any solution~$u$ of~(\ref{eq}) or of similar problems in the half-space $\Omega'=(0,+\infty)\times\R^{N-1}$ is well-defined and unique (roughly speaking, no oscillation occur). In some cases, we will prove that all solutions~$u$ converge as $x_1\to+\infty$ to the same limiting one-dimensional profile, irrespectively of~$u_0$. To do so, we will use two different approaches. The first one is a pure PDE approach based on comparisons with suitable sub-solutions and on Liouville type results. This paper indeed contains new Liouville type results of independent interest for the solutions of some elliptic equations in half-spaces $\R^{N-1}\times(0,+\infty)$ with homogeneous Dirichlet boundary conditions, or in the whole space~$\R^N$ (see Section~\ref{sec2} for more details). The second approach is a dynamical systems' approach which says that~$x_1$ can all the same be viewed as a time variable for a suitably defined dynamical system whose global attractor can be proved to exist and can be characterized. 

Let us now describe more precisely the types of assumptions we make on the functions~$f$, which are always assumed to be locally Lipschitz-continuous from $\R_+$ to~$\R$. The first class of functions we consider corresponds to functions~$f$ such that
\be\label{hypf1}\left\{\baa{l}
\exists\,\mu>0,\quad f>0\hbox{ on }(0,\mu),\quad f\le0\hbox{ on }[\mu,+\infty),\vspace{3pt}\\
\exists\,0<\mu'<\mu,\quad f\hbox{ is nonincreasing on }[\mu',\mu],\vspace{3pt}\\
\hbox{either }\Big[f(0)>0\Big]\hbox{ or }\Big[f(0)=0\hbox{ and }\displaystyle{\mathop{\liminf}_{s\to0^+}}\,\displaystyle{\frac{f(s)}{s}}>0\Big].\eaa\right.
\ee
Under assumption (\ref{hypf1}) on $f$, it is immediate to see that there exists a unique solution $V\in C^2(\R_+)$ of the one-dimensional equation
\be\label{eq1d}\left\{\baa{l}
V''(\xi)+f(V(\xi))=0\hbox{ for all }\xi\ge 0,\vspace{3pt}\\
V(0)=0<V(\xi)<\mu=V(+\infty)\hbox{ for all }\xi>0.\eaa\right.
\ee
Furthermore, $V'(\xi)>0$ for all $\xi\ge 0$.\par
Under assumption~(\ref{hypf1}), the behavior of the nontrivial solutions~$u$ of~(\ref{eq}) as $x_1\to+\infty$ is uniquely determined, as the following theorem shows.

\begin{theo}\label{th1}
Let $N$ be any integer such that $N\ge 2$ and assume that~$f$ satisfies~$(\ref{hypf1})$. Let~$u$ be any nonnegative and bounded solution of~$(\ref{eq})$, where $u_0:\R^{N-2}\times(0,+\infty)\to\R_+$ is any continuous and bounded function such that $u_0\not\equiv 0$ in $\R^{N-2}\times(0,+\infty)$. Then
\be\label{liminf}
\liminf_{R\to+\infty}\ \inf_{(R,+\infty)\times\R^{N-2}\times(R,+\infty)}u\ge\mu
\ee
and
\be\label{convV}
u(x_1+h,x',x_N)\to V(x_N)\ \hbox{ as }h\to+\infty\hbox{ in }C^{2,\beta}_b([A,+\infty)\times\R^{N-2}\times[0,B])
\ee
for all $A\in\R$, $B>0$ and $\beta\in[0,1)$, where $V\in C^2(\R_+)$ is the unique solution of $(\ref{eq1d})$.
\end{theo}

Notice that property (\ref{liminf}) means that the non-trivial nonnegative solutions~$u$ of~(\ref{eq}) are separated from~$0$, irrespectively of~$u_0$, far away from the boundary~$\partial\Omega$. If $u_0\le\mu$, then since $f\le 0$ on $[\mu,+\infty)$ and $u=\mu$ on $\partial\Omega^+$, where $\Omega^+=\Omega\cap\big\{u>\mu\}$, it follows that from the maximum principle applied in $\Omega^+$ (see~\cite{bcn2}, since $\R^N\backslash\overline{\Omega^+}$ contains the closure of an infinite open connected cone), that actually $u\le\mu$ in $\Omega^+$ whence $\Omega^+=\emptyset$ and $u\le\mu$ in~$\Omega$. In this case, it also follows from Theorem~\ref{th1} and standard elliptic estimates that the convergence~(\ref{convV}) holds not only locally in $x_N$, but in $C^{2,\beta}_b([A,+\infty)\times\R^{N-2}\times\R_+)$ for all $A\in\R$ and $\beta\in[0,1)$. However, without the assumption $u_0\le\mu$, it is not clear that this last convergence property holds globally with respect to $x_N$ in general.

The second class of functions $f$ we consider corresponds to the following assumption:
\be\label{hypf2}\left\{\baa{l}
f\ge 0\hbox{ on }\R_+,\vspace{3pt}\\
\forall\, z\in E,\ \displaystyle{\mathop{\liminf}_{s\to z^+}}\ \displaystyle{\frac{f(s)}{s-z}}>0,\eaa\right.
\ee
where
\be\label{defE}
E=\big\{z\in\R_+;\ f(z)=0\big\}
\ee
denotes the set of zeroes of $f$. A typical example of such a function~$f$ is $f(s)=|\sin s|$ for all $s\ge 0$, with $E=\pi\N$. More generally speaking, under the assumption~(\ref{hypf2}), it follows immediately that the set~$E$ is at most countable. Furthermore, it is easy to check that, for each $z\in E\backslash\{0\}$, there exists a unique solution $V_z\in C^2(\R_+)$ of the one-dimensional equation
\be\label{eq1dbis}\left\{\baa{l}
V_z''(\xi)+f(V_z(\xi))=0\hbox{ for all }\xi\ge 0,\vspace{3pt}\\
V_z(0)=0<V_z(\xi)<z=V_z(+\infty)\hbox{ for all }\xi>0.\eaa\right.
\ee
Furthermore, $V_z'(\xi)>0$ for all $\xi\ge 0$.\par
The following theorem states any solution of~(\ref{eq}) is asymptotically one-dimensional as $x_1\to+\infty$.

\begin{theo}\label{th2}
Let $N$ be any integer such that $N\ge 2$ and assume that~$f$
satisfies~$(\ref{hypf2})$. Let~$u$ be any nonnegative and bounded
solution of~$(\ref{eq})$, where
$u_0:\R^{N-2}\times(0,+\infty)\to\R_+$ is any continuous and
bounded function such that $u_0\not\equiv 0$ in
$\R^{N-2}\times(0,+\infty)$. Then there exists $R>0$ such that
$$\inf_{(R,+\infty)\times\R^{N-2}\times(R,+\infty)}u>0$$
and there exists $z\in E\backslash\{0\}$ such that
\be\label{uasympt}
u(x_1+h,x',x_N)\to V_z(x_N)\ \hbox{ as }h\to+\infty\hbox{ in }C^{2,\beta}_b([A,+\infty)\times\R^{N-2}\times\R_+)
\ee
for all $A\in\R$ and $\beta\in[0,1)$, where $V_z\in C^2(\R_+)$ is the unique solution of $(\ref{eq1dbis})$ with the limit $V_z(+\infty)=z$.
\end{theo}

This result shows that any non-trivial bounded solution~$u$ of~(\ref{eq}) converges to a single one-dimensional profile as~$x_1\to+\infty$. More precisely, given~$u$, the real number~$z$ defined by~(\ref{uasympt}) is unique and, in the proof of Theorem~\ref{th2}, the explicit expression of~$z$ will be provided. Observe that the asymptotic profile may now depend on the solution~$u$ (unlike in Theorem~\ref{th1}) but Theorem~\ref{th2} says that the oscillations in the~$x_1$~variable are excluded at infinity, for any solution~$u$.

The last two results are concerned with the analysis of the asymptotic behavior, as $x_1\to+\infty$, of the nonnegative bounded classical solutions~$u$ of
\be\label{eqbis}\left\{\baa{rcll}
\Delta u+f(u) & = & 0 & \hbox{ in }\Omega'=(0,+\infty)\times\R^{N-1},\vspace{3pt}\\
u(0,x_2,\ldots,x_N) & = & u_0(x_2,\ldots,x_N) & \hbox{ for all }(x_2,\ldots,x_N)\in\R^{N-1},\eaa\right.
\ee
in the half-space $\Omega'$, where the function $u_0:\R^{N-1}\to\R_+$ is given, continuous and bounded. The solutions~$u$ of~(\ref{eqbis}) are understood to be bounded, of class~$C^2(\Omega')$ and to be continuous on~$\overline{\Omega'}$. They are then automatically of class $C^{2,\beta}_b([\epsilon,+\infty)\times\R^{N-1})$ for all $\epsilon>0$ and $\beta\in[0,1)$. Firstly, under the same assumptions~(\ref{hypf2}) as in the previous theorem, the behavior as $x_1\to+\infty$ of any non-trivial solution~$u$ of~(\ref{eqbis}) is well-defined:

\begin{theo}\label{th3}
Let $N$ be any integer such that $N\ge 2$ and assume that~$f$ satisfies~$(\ref{hypf2})$. Let~$u$ be any nonnegative and bounded solution of~$(\ref{eqbis})$, where $u_0:\R^{N-1}\to\R_+$ is any continuous and bounded function such that $u_0\not\equiv 0$ in $\R^{N-1}$. Then there exists $z\in E\backslash\{0\}$ such that
\be\label{usasymptbis}
u(x_1+h,x_2,\ldots,x_N)\to z\ \hbox{ as }h\to+\infty\hbox{ in }C^{2,\beta}_b([A,+\infty)\times\R^{N-1})
\ee
for all $A\in\R$ and $\beta\in[0,1)$.
\end{theo}

Notice that the conclusion implies in particular that~$u$ is separated from~$0$ far away from the boundary $\{0\}\times\R^{N-1}$ of~$\Omega'$. Furthermore, as in Theorem~\ref{th2}, the real number~$z$ in~(\ref{usasymptbis}) is uniquely determined by~$u$ and its explicit value will be given during the proof. In Theorem~\ref{th3}, if instead of~(\ref{hypf2}) the function~$f$ now satisfies assumption~(\ref{hypf1}), then~$u$ may not converge in general to a constant as $x_1\to+\infty$. Furthermore, even if~$u$ does converge to a constant as $x_1\to+\infty$, that constant may not be equal to the real number~$\mu$ given in~(\ref{hypf1}). For instance, if there exists $\rho\in(\mu,+\infty)$ such that $f(\rho)=0$, then the constant function $u=\rho$ solves~(\ref{eq}) with $u_0=\rho$. Therefore, under assumption~(\ref{hypf1}), the asymptotic profile of a solution~$u$ of problem~(\ref{eqbis}) in the half-space~$\Omega'$ depends on~$u$ and is even not clearly well-defined in general. The situation is thus very different from Theorem~\ref{th1} about the existence and uniqueness of the asymptotic behavior of the solutions of problem~(\ref{eq}) in the quarter-space~$\Omega$.

However, under (\ref{hypf1}) and an additional appropriate assumption on $f$, the following result holds:

\begin{theo}\label{th4}
Let $N$ be any integer such that $N\ge 2$ and assume that, in addition to~$(\ref{hypf1})$,~$f$ is such that
\be\label{hypf3}
\liminf_{s\to z^-}\frac{f(s)}{s-z}>0\ \hbox{ for all }z>\mu\hbox{ such that }f(z)=0.
\ee
Let~$u$ be any nonnegative and bounded solution of~$(\ref{eqbis})$, where $u_0:\R^{N-1}\to\R_+$ is any continuous and bounded function such that $u_0\not\equiv 0$ in $\R^{N-1}$. Then there exists $z\ge\mu$ such that $f(z)=0$ and
$$u(x_1+h,x_2,\ldots,x_N)\to z\ \hbox{ as }h\to+\infty\hbox{ in }C^{2,\beta}_b([A,+\infty)\times\R^{N-1})$$
for all $A\in\R$ and $\beta\in[0,1)$.
\end{theo}

\begin{rem}{\rm It is worth noticing that all above results hold in any dimension $N\ge 2$.}
\end{rem}


\section{The PDE approach}\label{sec2}

In this section, we use a pure PDE approach to prove the main results announced in Section~\ref{intro}. In Section~\ref{sec21}, we deal with the case of problem~(\ref{eq}) set in the quarter-space~$\Omega=(0,+\infty)\times\R^{N-2}\times(0,+\infty)$, while Section~\ref{sec22} is concerned with problem~(\ref{eqbis}) set in the half-space~$\Omega'=(0,+\infty)\times\R^{N-1}$.


\subsection{Problem~(\ref{eq}) in the quarter-space~$\Omega\!=\!(0,+\infty)\!\times\!\R^{N-2}\!\times\!(0,+\infty)$}\label{sec21}

Let us first begin with the\hfill\break

\noindent{\bf{Proof of Theorem~\ref{th1}.}} The proof is divided into two main steps: we first prove that~$u$ is bounded from below away from $0$ when $x_1$ and $x_N$ are large, uniformly with respect to~$x'\in\R^{N-2}$. Then, we pass to the limit as $x_1\to+\infty$ and use a classification result, which leads to the uniqueness and one-dimensional symmetry of the limiting profiles of $u$ as $x_1\to+\infty$.\par
First of all, observe that, since $f(0)\ge 0$, the strong maximum principle implies that the function~$u$ is either positive in $\Omega$, or identically equal to~$0$ in~$\Omega$. But since $u$ is continuous up to $\{0\}\times\R^{N-2}\times(0,+\infty)$ and since $u_0\not\equiv 0$, it follows that $u>0$ in $\Omega$.\par
{\it{Step 1.}} In the sequel, for any $x\in\R^N$ and $R>0$, denote $B(x,R)$ the open euclidean ball of centre $x$ and radius $R$. For each $R>0$, let $\lambda_R$ be the principal eigenvalue of the Laplace operator in $B(0,R)$ with Dirichlet boundary condition on $\partial B(0,R)$, and let $\varphi_R$ be the normalized principal eigenfunction, that is
\be\label{phiR}\left\{\baa{ll}
\Delta\varphi_R+\lambda_R\,\varphi_R=0 & \hbox{ in }B(0,R),\vspace{3pt}\\
\varphi_R>0 & \hbox{ in }B(0,R),\vspace{3pt}\\
\|\varphi_R\|_{L^{\infty}(B(0,R))}=\varphi_R(0)=1, & \vspace{3pt}\\
\varphi_R=0 & \hbox{ on }\partial B(0,R).\eaa\right.
\ee\par
Notice that $\lambda_R\to 0$ as $R\to+\infty$. If $f(0)=0$, we can then choose $R>0$ large enough so that
$$\lambda_R<\liminf_{s\to 0^+}\frac{f(s)}{s}.$$
If $f(0)>0$, we simply choose $R=1$. Then, fix a point $x_0\in\Omega$ in such a way that~$\overline{B(x_0,R)}\subset\Omega$. Since $u$ is continuous and positive on $\overline{B(x_0,R)}$, there holds~$\min_{\overline{B(x_0,R)}}u>0$. Therefore, it follows from the choice of $R$ that there exists $\epsilon>0$ small enough, such that the function
$$\underline{u}(x)=\epsilon\,\varphi_R(x-x_0)$$
is a subsolution in $\overline{B(x_0,R)}$, that is
\be\label{subsolution}
\Delta\underline{u}+f(\underline{u})\ge 0\ \hbox{ and }\ \underline{u}<u\ \hbox{ in }\ \overline{B(x_0,R)}.
\ee\par
Next, let $\widetilde{x}_0$ be any point in $\Omega$ such that
$$\overline{B(\widetilde{x}_0,R)}\subset\overline{\Omega},$$
that is $\widetilde{x}_0=(\widetilde{x}_{0,1},\widetilde{x}'_0,\widetilde{x}_{0,N})$ with $\widetilde{x}_{0,1}\ge R$ and $\widetilde{x}_{0,N}\ge R$. For all $t\in[0,1]$, call
$$y_t=x_0+t\,(\widetilde{x}_0-x_0)$$
and observe that $\overline{B(y_t,R)}\subset\Omega$ for all $t\in[0,1)$. Define
$$\underline{u}_t(x)=\underline{u}(x-y_t+x_0)=\epsilon\,\varphi_R(x-y_t)$$
for all $t\in[0,1]$ and $x\in\overline{B(y_t,R)}$. By continuity and from (\ref{subsolution}), there holds $\underline{u}_t<u$ in~$\overline{B(y_t,R)}$ for $t\in[0,t_0]$, where $t_0>0$ is small enough. On the other hand, for all $t\in[0,1]$, the function $\underline{u}_t$ is a subsolution of the equation satisfied by $u$, that is
$$\Delta\underline{u}_t+f(\underline{u}_t)\ge 0\hbox{ in }\overline{B(y_t,R)}.$$
We shall now use a sliding method (see \cite{bcn2,bn}) to conclude that
$$\underline{u}_t<u\hbox{ in }\overline{B(y_t,R)}\hbox{ for all }t\in[0,1).$$
Indeed, if this were not true, there would then exist a real number $t^*\in(0,1)$ such that the inequality $\underline{u}_{t^*}\le u$ holds in $\overline{B(y_{t^*},R)}$ with equality at some point $x^*\in\overline{B(y_{t^*},R)}$. Since~$\overline{B(y_{t^*},R)}\subset\Omega$, $u>0$ in $\Omega$ and $\underline{u}_{t^*}=0$ on $\partial B(y_{t^*},R)$, one has $x^*\in B(y_{t^*},R)$. But since~$\underline{u}_{t^*}$ is a subsolution of the equation satisfied by $u$, the strong maximum principle yields~$\underline{u}_{t^*}=u$ in $B(y_{t^*},R)$ and also on the boundary by continuity, which is impossible. One has then reached a contradiction. Hence, $\underline{u}_t<u$ in $\overline{B(y_t,R)}$ for all $t\in[0,1)$. By continuity, one also gets that $\underline{u}_1\le u$ in $\overline{B(y_1,R)}$.\par
Eventually, for any $\widetilde{x}_0=(\widetilde{x}_{0,1},\widetilde{x}'_0,\widetilde{x}_{0,N})$ with $\widetilde{x}_{0,1}\ge R$ and $\widetilde{x}_{0,N}\ge R$, there holds
$$u(\widetilde{x}_0)\ge\underline{u}_1(\widetilde{x}_0)=\epsilon\,\varphi_R(0)=\epsilon.$$
In other words,
\be\label{ineqeps}
u\ge\epsilon\ \hbox{ in }[R,+\infty)\times\R^{N-2}\times[R,+\infty).\ee\par
{\it{Step 2.}} Let $(x_{1,n})_{n\in\N}$ be any sequence of positive numbers such that $x_{1,n}\to+\infty$ as~$n\to+\infty$, and let $(x'_n)_{n\in\N}$ be any sequence in $\R^{N-2}$. From standard elliptic estimates, there exists a subsequence such that the functions
$$u_n(x)=u(x_1+x_{1,n},x'+x'_n,x_N)$$
converge in $C^{2,\beta}_{loc}(\R^N_+)$, for all $\beta\in[0,1)$, to a bounded classical solution $u_{\infty}$ of
\be\label{eqlim}\left\{\baa{rcll}
\Delta u_{\infty}+f(u_{\infty}) & = & 0 & \hbox{ in }\R^N_+,\vspace{3pt}\\
u_{\infty} & = & 0 & \hbox{ on }\partial\R^N_+,\eaa\right.
\ee
where $\R^N_+=\R^{N-1}\times[0,+\infty)$. Furthermore, $u_{\infty}\ge 0$ and $u_{\infty}\not\equiv 0$ in $\R^N_+$ from Step~1, since
$$u_{\infty}\ge\epsilon>0\ \hbox{ in }\ \R^{N-1}\times[R,+\infty)$$
from~(\ref{ineqeps}). Thus, $u_{\infty}>0$ in $\R^{N-1}\times(0,+\infty)$ from the strong maximum principle. It follows from Theorems~1.1 and~1.2 of Berestycki, Caffarelli and Nirenberg~\cite{bcn2}\footnote{In~\cite{bcn2}, the function $f$ was assumed to be globally Lipschitz-continuous. Here,~$f$ is just assumed to be locally Lipschitz-continuous. However, since $u$ is bounded, it is always possible to find a Lipschitz-continuous function~$\tilde{f}:\R^+\to\R$ satisfying~(\ref{hypf1}) and such that~$\tilde{f}$ and~$f$ coincide on the range of~$u$.} (see also~\cite{a,cs}) that $u_{\infty}$ is unique and has one-dimensional symmetry. By uniqueness of the problem (\ref{eq1d}), one gets that $u_{\infty}(x)=V(x_N)$ for all $x\in\R^N_+$, and the limit does not depend on the sequences $(x_{1,n})_{n\in\N}$ or $(x'_n)_{n\in\N}$. Property~(\ref{convV}) of Theorem~\ref{th1} then follows from the uniqueness of the limit.\par
{\it{Step 3.}} Let us now prove formula~(\ref{liminf}). One already knows from~(\ref{ineqeps}) that
$$m:=\liminf_{R\to+\infty}\ \inf_{(R,+\infty)\times\R^{N-2}\times(R,+\infty)}u\ge\epsilon>0.$$
Let $(x_n)_{n\in\N}=(x_{1,n},x'_n,x_{N,n})_{n\in\N}$ be a sequence in $\Omega$ such that $(x_{1,n},x_{N,n})_{n\in\N}\to(+\infty,+\infty)$ and $u(x_n)\to m$ as $n\to+\infty$. Up to extraction of a subsequence, the functions
$$v_n(x)=u(x+x_n)$$
converge in $C^2_{loc}(\R^N)$ to a classical bounded solution $v_{\infty}$ of
$$\Delta v_{\infty}+f(v_{\infty})=0\ \hbox{ in }\R^N$$
such that $v_{\infty}\ge m$ in $\R^N$ and $v_{\infty}(0)=m>0$. Thus, $f(m)\le 0$, whence $m\ge\mu$ due to~(\ref{hypf1}). The proof of Theorem~\ref{th1} is thereby complete.\hfill$\Box$\break

The proof of Theorem~\ref{th2} is based on two Liouville type results for the bounded nonnegative solutions~$u$ of the elliptic equation
$$\Delta u+f(u)=0$$
in the whole space $\R^N$ or in the half-space $\R^{N-1}\times\R_+$ with Dirichlet boundary condition on $\R^{N-1}\times\{0\}$.

\begin{theo}\label{thspace}
Let $N$ be any integer such that $N\ge 1$ and assume that the function~$f$ satisfies~$(\ref{hypf2})$. Let~$u$ be a bounded nonnegative solution of
\be\label{eqRN}
\Delta u+f(u)=0\ \hbox{ in }\R^N.
\ee
Then $u$ is constant.
\end{theo}

The following result is concerned with the one-dimensional symmetry of nonnegative bounded solutions in a half-space with Dirichlet boundary conditions.

\begin{theo}\label{thhalfspace}
Let $N$ be any integer such that $N\ge 1$ and assume that the function~$f$ satisfies~$(\ref{hypf2})$. Let~$u$ be a bounded nonnegative solution of
\be\label{eqhalf}\left\{\baa{rcl}
\Delta u+f(u) & = & 0\ \hbox{ in }\R^N_+=\R^{N-1}\times\R_+,\vspace{3pt}\\
u & = & 0\ \hbox{ on }\partial\R^N_+=\R^{N-1}\times\{0\}.\eaa\right.
\ee
Then $u$ is a function of $x_N$ only. Furthermore, either $u=0$ in $\R^{N-1}\times\R_+$ or there exists~$z>0$ such that $f(z)=0$ and $u(x)=V_z(x_N)$ for all $x\in\R^{N-1}\times\R_+$, where the function~$V_z$ satisfies equation~$(\ref{eq1dbis})$.
\end{theo}

These results are of independent interest and will be proved in Section~\ref{secliouville}. Notice that one of the main points is that they hold in any dimensions~$N\ge 1$ without any other assumption on~$u$ than its boundedness. In low dimensions $N\le 4$, and under the additional assumption that $u$ is stable, the conclusion of Theorem~\ref{thspace} holds for any nonnegative function~$f$ of class~$C^1(\R_+)$, see Dupaigne and Farina~\cite{df}. Consequently, because of the monotonicity result in the direction~$x_N$ due to Berestycki, Caffarelli and Nirenberg~\cite{bcn2} and Dancer \cite{d} (since $f(0)\ge 0$), it follows that the conclusion of Theorem~\ref{thhalfspace} holds for any nonnegative function~$f$ of class~$C^1(\R_+)$, provided that $N\le 5$, see Farina and Valdinoci~\cite{fv}. However, observe that the nonnegativity and the $C^1$ character of $f$ are incompatible with~(\ref{hypf2}) for any posi\-tive zero~$z$ of~$f$. Furthermore, assumption~(\ref{hypf2}) is crucially used in the proof of Theorem~\ref{thspace} and~\ref{thhalfspace}. It is actually not true that these theorems stay valid in general when~$f$ is just assumed to be nonnegative and locally Lipschitz-continuous. For instance, non-constant solutions of~(\ref{eqRN}), which are even stable, exist for power-like nonlinearities~$f$ in high dimensions (see~\cite{f} and the references therein).\hfill\break\par
With these results in hand, let us turn to the\hfill\break

\noindent{\bf{Proof of Theorem~\ref{th2}.}} Observe that, from~(\ref{hypf2}), either $f(0)>0$, or $f(0)=0$ and $\liminf_{s\to 0^+}f(s)/s>0$. Therefore, as in the proof of Theorem~\ref{th1}, there exist~$R>0$ and~$\epsilon>0$ such that
\be\label{uepsilon}
u\ge\epsilon\ \hbox{ in }[R,+\infty)\times\R^{N-2}\times[R,+\infty).
\ee
Set
\be\label{defM}
M=\lim_{A\to+\infty}\,\sup_{[A,+\infty)\times\R^{N-2}\times[0,+\infty)}\,u.
\ee
Our goal is to prove that the conclusion of Theorem~\ref{th2} holds with $z=M$.\par
Since $u$ is bounded and satisfies~(\ref{uepsilon}),~$M$ is such that $\epsilon\le M<+\infty$. Furthermore, there exists a sequence~$(x_n)_{n\in\N}=(x_{1,n},x'_n,x_{N,n})_{n\in\N}$ of points in $\overline{\Omega}$ such that $x_{1,n}\to+\infty$ and $u(x_n)\to M$ as~$n\to+\infty$.\par
Assume first, up to extraction of a subsequence, that the sequence $(x_{N,n})_{n\in\N}$ converges to a, nonnegative, real number $x_{N,\infty}$ as $n\to+\infty$. From standard elliptic estimates, the functions
$$u_n(x)=u(x_1+x_{1,n},x'+x'_n,x_N)$$
converge in $C^{2,\beta}_{loc}(\R^{N-1}\times\R_+)$ for all $\beta\in[0,1)$, up to extraction of another subsequence, to a bounded nonnegative solution $u_{\infty}$ of the problem~(\ref{eqhalf}) in the half-space $\R^{N-1}\times[0,+\infty)$, such that
$$u_{\infty}(0,0,x_{N,\infty})=M=\sup_{\R^{N-1}\times[0,+\infty)}\,u_{\infty}>0.$$
It follows from Theorem~\ref{thhalfspace} that $u_{\infty}(x)=V_z(x_N)$ is a one-dimensional increasing solution of~(\ref{eq1dbis}), whence $z=M$. But since $V_M$ is (strictly) increasing, it cannot reach its maximum~$M$ at the finite point $x_{N,\infty}$. This case is then impossible.\par
Therefore, one can assume without loss of generality that $x_{N,n}\to+\infty$ as $n\to+\infty$. From standard elliptic estimates, the functions
$$u_n(x)=u(x+x_n)$$
converge in $C^{2,\beta}_{loc}(\R^N)$ for all $\beta\in[0,1)$, up to extraction of another subsequence, to a bounded nonegative solution $u_{\infty}$ of the problem~(\ref{eqRN}) in $\R^N$, such that
$$u_{\infty}(0)=M=\sup_{\R^N}\,u_{\infty}>0.$$
Theorem~\ref{thspace} implies that $u_{\infty}=M$ in $\R^N$, whence $f(M)=0$. Furthermore, since the limit~$M$ is unique, the convergence of the functions~$u_n$ to the constant~$M$ holds for the whole sequence.\par
Now, in order to complete the proof of Theorem~\ref{th2}, we shall make use of the following lemma of independent interest:

\begin{lem}\label{lemgv} Let $g:\R_+\to\R$ be a locally Lipschitz-continuous nonnegative function. Then, for each $z>0$ such that $g(z)=0$ and for each $\epsilon\in(0,z]$, there exist $R'=R'_{g,z,\epsilon}>0$ and a classical solution~$v$ of
\be\label{eqv}\left\{\baa{rcl}
\Delta v+g(v) &= & 0\hbox{ in }\overline{B(0,R')},\vspace{3pt}\\
0\ \le\ v & < & z\hbox{ in }\overline{B(0,R')},\vspace{3pt}\\
v & = & 0\hbox{ on }\partial B(0,R'),\\
v(0)\,=\displaystyle{\mathop{\max}_{\overline{B(0,R')}}}\,v & \ge & z-\epsilon.\eaa\right.
\ee
\end{lem}

The proof of this lemma is postponed at the end of this section. Let us now finish the proof of Theorem~\ref{th2}. Fix an arbitrary~$\epsilon$ in~$(0,M]$. Let $R'(\epsilon)=R'_{f,M,\epsilon}$ be as in Lemma~\ref{lemgv} and let $v$ be a solution of~(\ref{eqv}) with $g=f$ and $z=M$. Since the functions $u_n(x)=u(x+x_n)$ converge locally uniformly in $\R^N$ to the constant $M$ and since
$$\max_{\overline{B(0,R'(\epsilon))}}v<M,$$
there exists $n_0\in\N$ large enough so that $\overline{B(x_{n_0},R'(\epsilon))}\subset\Omega$ and
\be\label{ineqvu}
v(x-x_{n_0})<u(x)\hbox{ for all }x\in\overline{B(x_{n_0},R'(\epsilon))}.
\ee
But since $v$ solves the same elliptic equation as~$u$, the same sliding method as in Theorem~\ref{th1} implies that
\be\label{uvepsilon}
u\ge v(0)\ge M-\epsilon\ \hbox{ in }[R'(\epsilon),+\infty)\times\R^{N-2}\times[R'(\epsilon),+\infty).
\ee\par
Lastly, choose any sequence $(\widetilde{x}_{1,n})_{n\in\N}$ converging to $+\infty$, and any sequence~$(\widetilde{x}'_n)_{n\in\N}$ in~$\R^{N-2}$. Up to extraction of a subsequence, the functions
$$\widetilde{u}_n(x)=u(x_1+\widetilde{x}_{1,n},x'+\widetilde{x}'_n,x_N)$$
converge in $C^{2,\beta}_{loc}(\R^{N-1}\times\R_+)$ for all $\beta\in[0,1)$ to a bounded nonnegative solution~$\widetilde{u}_{\infty}$ of problem~(\ref{eqhalf}) in the half-space $\R^{N-1}\times\R_+$. The function $\widetilde{u}_{\infty}$ satisfies $\widetilde{u}_{\infty}\le M$ in $\R^{N-1}\times\R_+$ by definition of~$M$, while
$$\lim_{A\to+\infty}\ \inf_{\R^{N-1}\times[A,+\infty)}\,\widetilde{u}_{\infty}\ge M$$
because $\epsilon>0$ in~(\ref{uvepsilon}) can be arbitrarily small. Theorem~\ref{thhalfspace} and the above estimates imply that
$$\widetilde{u}_{\infty}(x)=V_M(x_N)\ \hbox{ for all }x\in\R^{N-1}\times\R_+.$$
Since this limit does not depend on any subsequence, and due to~(\ref{defM}),~(\ref{uvepsilon}) and standard elliptic estimates, it follows in particular that
$$u(x_1+h,x',x_N)\to V_M(x_N)\ \hbox{ in }C^{2,\beta}_b([A,+\infty)\times\R^{N-2}\times\R_+)\ \hbox{ as }h\to+\infty,$$
for all $A\in\R$ and $\beta\in[0,1)$. The proof of Theorem~\ref{th2} is thereby complete.\hfill$\Box$

\begin{rem}{\rm Instead of Theorem~\ref{thspace}, if $f$ is just assumed to be nonnegative and locally Lipschitz-continuous on $\R_+$ and if~$u_{\infty}$ is a solution of~(\ref{eqRN}) which reaches its maximum and is such that $f(\max_{\R^N}u_{\infty})=0$ (as in some assumptions of~\cite{bcn1,bcn3}), then~$u_{\infty}$ is constant, from the strong maximum principle. Therefore, it follows from similar arguments as in the above proof that if, instead of~(\ref{hypf2}) in Theorem~\ref{th2}, the function~$f$ is just assumed to be nonnegative, locally Lipschitz-continuous and positive almost everywhere on~$\R_+$ and if~$u$ is a bounded nonnegative solution of~(\ref{eq}) such that~$f(M)=0$, where~$M$ is defined by~(\ref{defM}), then either $M=0$ and $u(x_1+h,x',x_N)\to 0$ in $C^2([A,+\infty)\times\R^{N-2}\times\R_+)$ as $h\to+\infty$ for all $A\in\R$, or the conclusion~(\ref{uasympt}) holds with~$z=M$.}
\end{rem}

\noindent{\bf{Proof of Lemma~\ref{lemgv}.}} The proof uses classical variational arguments, which we sketch here for the sake of completeness (see also e.g.~\cite{bhm1,eim} for applications of this method). Let~$g$ and~$z$ be as in Lemma~\ref{lemgv} and let $\widetilde{g}$ be the function defined in $\R$ by
$$\tilde{g}(s)=\left\{\baa{ll}
g(0) & \hbox{ if }s<0,\vspace{3pt}\\
g(s) & \hbox{ if }0\le s\le z,\vspace{3pt}\\
0 & \hbox{ if }s>z.\eaa\right.$$
The function $\widetilde{g}$ is nonnegative, bounded and Lipschitz-continuous on $\R$. Set
$$G(s)=\int_s^z\widetilde{g}(\tau)\,d\tau\ge 0$$
for all $s\in\R$. The function $G$ is nonnegative and Lipschitz-continuous on $\R$.\par
Let $r$ be any positive real number. Define
$$I_{r}(v)=\frac{1}{2}\int_{B(0,r)}|\nabla v|^2+\int_{B(0,r)}G(v)$$
for all $v\in H^1_0(B(0,r))$. The functional $I_{r}$ is well-defined in $H^1_0(B(0,r))$ and it is coercive, from Poincar\'e's inequality and the nonnegativity of~$G$. From Rellich's and Lebesgue's theorems, the functional~$I_{r}$ has a minimum~$v_{r}$ in~$H^1_0(B(0,r))$. The function~$v_{r}$ is a weak and hence, from the elliptic regularity theory, a classical $C^2(\overline{B(0,r)})$ solution of the equation
$$\left\{\baa{rcl}
\Delta v_{r}+\widetilde{g}(v_{r}) &= & 0\hbox{ in }\overline{B(0,r)},\vspace{3pt}\\
v_{r} & = & 0\hbox{ on }\partial B(0,r).\eaa.\right.$$
Since $\widetilde{g}\ge 0$ on $(-\infty,0]$, it follows from the strong maximum principle that $v_{r}\ge 0$ in~$\overline{B(0,r)}$. Furthermore, either $v_{r}=0$ in $\overline{B(0,r)}$, or $v_{r}>0$ in $B(0,r)$. Similarly, since $\widetilde{g}=0$ on~$[z,+\infty)$, one gets that $v_{r}<z$ in $\overline{B(0,r)}$. Consequently, $\widetilde{g}(v_{r})=g(v_{r})$ in~$\overline{B(0,r)}$. It also follows from the method of moving planes and Gidas, Ni and Nirenberg~\cite{gnn} that $v_{r}$ is radially symmetric and decreasing with respect to $|x|$ (provided that $v_{r}\not\equiv 0$ in~$\overline{B(0,r)}$). In all cases, there holds
$$0\le v_{r}(0)=\max_{\overline{B(0,r)}}v_{r}<z.$$\par
In order to complete the proof of Lemma~\ref{lemgv}, it is sufficient to prove that, given~$\epsilon$ in~$(0,z]$, there exists $r>0$ such that $v_r(0)\ge z-\epsilon$. Let $\epsilon\in(0,z]$ and assume that~$\max_{\overline{B(0,r)}}v_{r}=v_{r}(0)<z-\epsilon$ for all $r>0$. Observe that the function $G$ is nonincreasing in~$\R$, and actually decreasing and positive on the interval $[0,z)$, from~(\ref{hypf2}). Therefore,
\be\label{Irv}
I_r(v_r)\ge\alpha_N\,r^N\,G(z-\epsilon)
\ee
for all $r>0$, where $\alpha_N>0$ denotes the Lebesgue measure of the unit euclidean ball in~$\R^N$. For $r>1$, let $w_r$ be the test function defined in $\overline{B(0,r)}$ by
$$w_r(x)=\left\{\baa{ll}
z & \hbox{ if }|x|<r-1,\vspace{3pt}\\
z\,(r-|x|) & \hbox{ if }r-1\le|x|\le r.\eaa\right.$$
This function $w_r$ belongs to $H^1_0(B(0,r))$ and $|\nabla w_r|^2$ and $G(w_r)$ are supported on the shell~$\overline{B(0,r)}\backslash B(0,r-1)$. Thus, there exists a constant $C$ independent of $r$ such that
\be\label{Irw}
I_r(w_r)\le C\,(r^n-(r-1)^n)
\ee
for all $r>1$. But since $I_r(v_r)\le I_r(w_r)$, by definition of $v_r$, and since $G(z-\epsilon)>0$, inequalities~(\ref{Irv}) and~(\ref{Irw}) lead to a contradiction as $r\to+\infty$. Therefore, there exists a radius $R'>0$ such that $v_{R'}(0)\ge z-\epsilon$, and $v_{R'}$ solves~(\ref{eqv}). The proof of Lemma~\ref{lemgv} is now complete.\hfill$\Box$


\subsection{Problem~(\ref{eqbis}) in the half-space~$\Omega'=(0,+\infty)\times\R^{N-1}$}\label{sec22}

Let us now turn to problem~(\ref{eqbis}) set in the half-space $\Omega'=(0,+\infty)\times\R^{N-1}$. This section is devoted to the proof of Theorems~\ref{th3} and~\ref{th4}.\hfill\break

\noindent{\bf{Proof of Theorem~\ref{th3}.}} Assume here that~$f$ satisfies~(\ref{hypf2}). First of all, as in Step~1 of the proof of Theorem~\ref{th1}, there exist $\epsilon>0$ and $R>0$ such that
\be\label{ineqhalf}
u\ge\epsilon\ \hbox{ in }[R,+\infty)\times\R^{N-1}.
\ee
Call
$$M=\lim_{A\to+\infty}\ \sup_{[A,+\infty)\times\R^{N-1}}u.$$
We shall now prove that the conclusion of Theorem~\ref{th3} holds with $z=M$. Choose a sequence $(x_n)_{n\in\N}=(x_{1,n},\ldots,x_{N,n})_{n\in\N}$ in~$\Omega'$ such that $x_{1,n}\to+\infty$ and $u(x_n)\to M$ as $n\to+\infty$. As in the proof of Theorem~\ref{th2}, it follows from Theorem~\ref{thspace} that the functions
$$u_n(x)=u(x+x_n)$$
converge in $C^{2,\beta}_{loc}(\R^N)$ for all $\beta\in[0,1)$ to the constant~$M$, whence $f(M)=0$.\par
Then, for any $\epsilon\in(0,M]$, let $R'(\epsilon)=R'_{f,M,\epsilon}$ be as in Lemma~\ref{lemgv} and let~$v$ be a solution of~(\ref{eqv}) with $g=f$ and $z=M$.  As in the proof of Theorem~\ref{th2}, there exists $n_0\in\N$ large enough so that $\overline{B(x_{n_0},R'(\epsilon))}\subset\Omega'$ and~(\ref{ineqvu}) holds. The sliding method yields
$$u\ge v(0)\ge M-\epsilon\ \hbox{ in }[R'(\epsilon),+\infty)\times\R^{N-1}.$$
Since $\epsilon>0$ is arbitrarily small, the definition of~$M$ implies that, for all $A\in\R$,
$$u(x_1+h,x_2,\ldots,x_N)\to M\ \hbox{ as }h\to+\infty$$
uniformly with respect to $(x_1,x_2,\ldots,x_N)\in[A,+\infty)\times\R^{N-1}$. The convergence also holds in $C^{2,\beta}_b([A,+\infty)\times\R^{N-1})$ for all $\beta\in[0,1)$ from standard elliptic estimates. The proof of Theorem~\ref{th3} is thereby complete.\hfill$\Box$\break

In the case when $f$ satisfies~(\ref{hypf1}) and~(\ref{hypf3}), the conclusion is similar to that of Theorem~\ref{th3}, as the following proof of Theorem~\ref{th4} will show. As a matter of fact, it is also based on a Liouville type result for the bounded nonnegative solutions of~(\ref{eqRN}), which is the counterpart of Theorem~\ref{thspace} under assumptions~(\ref{hypf1}) and~(\ref{hypf3}).

\begin{theo}\label{thspacebis}
Let $N$ be any integer such that $N\ge 1$ and assume that the function~$f$ satisfies~$(\ref{hypf1})$ and~$(\ref{hypf3})$. Then any bounded nonnegative solution~$u$ of~$(\ref{eqRN})$ is constant.
\end{theo}

The proof is postponed in Section~\ref{secliouville} and we now complete the\hfill\break

\noindent{\bf{Proof of Theorem~\ref{th4}.}} First of all, as in Step~1 of Theorem~\ref{th1}, there exist $\epsilon>0$ and $R>0$ such that~(\ref{ineqhalf}) holds. Call now
$$m=\lim_{A\to+\infty}\ \inf_{[A,+\infty)\times\R^{N-1}}u.$$
One has $m\in[\epsilon,+\infty)$. Let $(x_n)_{n\in\N}=(x_{1,n},\ldots,x_{N,n})_{n\in\N}$ be a sequence in~$\Omega'$ such that $x_{1,n}\to+\infty$ and $u(x_n)\to m$ as $n\to+\infty$. Up to extraction of a subsequence, the functions
$$u_n(x)=u(x+x_n)$$
converge in $C^{2,\beta}_{loc}(\R^N)$ for all $\beta\in[0,1)$ to a classical bounded solution~$u_{\infty}$ of~(\ref{eqRN}) such that $u_{\infty}\ge m>0$ in~$\R^N$ and $u_{\infty}(0)=m$. Theorem~\ref{thspacebis} then implies that~$u_{\infty}$ is constant in~$\R^N$, wence it is identically equal to~$m$ and $f(m)=0$.\par
Call now
$$M'=\sup_{\Omega'}u$$
and let $g:[0,+\infty)\to\R$ be the function defined by
$$g(s)=\left\{\baa{ll}
-f(M'+1-s) & \hbox{if }0\le s\le M'+1-m,\vspace{3pt}\\
0 & \hbox{if }s>M'+1-m.\eaa\right.$$
The function $g$ is Lipschitz-continuous and nonnegative. The real number
$$z=M'+1-m$$
is positive and fulfills $g(z)=f(m)=0$. Choose any $\epsilon$ in $(0,z]$. From Lemma~\ref{lemgv}, there exist $R'>0$ and a classical solution~$v$ of~(\ref{eqv}) in $\overline{B(0,R')}$, that is
$$\left\{\baa{rcl}
\Delta v+g(v) &= & 0\hbox{ in }\overline{B(0,R')},\vspace{3pt}\\
0\ \le\ v & < & z\hbox{ in }\overline{B(0,R')},\vspace{3pt}\\
v & = & 0\hbox{ on }\partial B(0,R'),\\
v(0)\,=\displaystyle{\mathop{\max}_{\overline{B(0,R')}}}\,v & \ge & z-\epsilon.\eaa\right.$$
The function $V=M'+1-v$ then satisfies
$$\left\{\baa{rcl}
\Delta V+f(V) &= & 0\hbox{ in }\overline{B(0,R')},\vspace{3pt}\\
m\ <\ V & \le & M'+1\hbox{ in }\overline{B(0,R')},\vspace{3pt}\\
V & = & M'+1\hbox{ on }\partial B(0,R'),\\
V(0)\,=\displaystyle{\mathop{\min}_{\overline{B(0,R')}}}\,V & \le & m+\epsilon.\eaa\right.$$
Since $u_n(x)=u(x+x_n)\to u_{\infty}(x)=m$ as $n\to+\infty$ in $C^{2,\beta}_{loc}(\R^N)$ for all $\beta\in[0,1)$, it follows that there exists $n_0\in\N$ large enough so that $\overline{B(x_{n_0},R')}\subset\Omega'$ and
$$V(x-x_{n_0})>u(x)\ \hbox{ for all }x\in\overline{B(x_{n_0},R')}.$$
Since $V=\sup_{\Omega'}u+1>\sup_{\Omega'}u$ on $\partial B(0,R')$, it follows from the elliptic maximum principle and the sliding method that
$$u(x)\le V(0)\le m+\epsilon\ \hbox{ for all }x\in[R',+\infty)\times\R^{N-1}.$$
Owing to the definition of~$m$, one concludes that, for all $A\in\R$,
$$u(x_1+h,x_2,\ldots,x_N)\to m\ \hbox{ as }h\to+\infty$$
uniformly with respect to $(x_1,x_2,\ldots,x_N)\in[A,+\infty)\times\R^{N-1}$ and the convergence holds in~$C^{2,\beta}_b([A,+\infty)\times\R^{N-1})$ for all $\beta\in[0,1)$ from standard elliptic estimates. The proof of Theorem~\ref{th4} is thereby complete.\hfill$\Box$


\section{Classification results in the whole space~$\R^N$ or in the half-space $\R^{N-1}\times(0,+\infty)$ with Dirichlet boundary conditions}\label{secliouville}

This section is devoted to the proof of the Liouville type results for the bounded nonnegative solutions~$u$ of problems~(\ref{eqRN}) or~(\ref{eqhalf}). Theorems~\ref{thhalfspace} and~\ref{thspacebis} are actually corollaries of Theorem~\ref{thspace}. We then begin with the proof of the latter.\hfill\break

\noindent{\bf{Proof of Theorem~\ref{thspace}.}} Let $u$ be a bounded nonnegative solution of~(\ref{eqRN}) under assumption~(\ref{hypf2}). Denote
$$m=\inf_{\R^N}\,u\ \ge 0.$$
Since $f(m)\ge 0$, the constant $m$ is a subsolution for~(\ref{eqRN}). It follows from the strong elliptic maximum principle that either $u=m$ in $\R^N$, or $u>m$ in $\R^N$.\par
Let us prove that the second case, that is $u>m$, is impossible. That will give the desired conclusion. Assume that $u>m$ in $\R^N$ and let us get a contradiction. Let us first check that
\be\label{mzero}
f(m)=0.
\ee
This could be done by considering a sequence along which~$u$ converges to its minimum; after changing the origin, the limiting function would be identically equal to $m$ from the strong maximum principle, which would yield~(\ref{mzero}). Let us choose an alternate elementary parabolic argument. Assume that $f(m)>0$ and let $\xi:[0,T)\to\R$ be the maximal solution of
$$\left\{\baa{rcl}
\xi'(t) & = & f(\xi(t))\ \hbox{ for all }t\in[0,T),\vspace{3pt}\\
\xi(0) & = & m.\eaa\right.$$
The maximal existence time~$T$ satisfies $0<T\le+\infty$ (and $T=+\infty$ if $f$ is globally Lipschitz-continuous). Since $\xi(0)\le u$ in $\R^N$, it follows from the parabolic maximum principle for the equation $v_t=\Delta v+f(v)$, satisfied by both $\xi$ and $u$ in $[0,T)\times\R^N$, that
$$\xi(t)\le u(x)\hbox{ for all }x\in\R^N\hbox{ and }t\in[0,T).$$
But $\xi'(0)=f(\xi(0))=f(m)>0$. Hence, there exists $\tau\in(0,T)$ such that $\xi(\tau)>m$, whence $u(x)\ge\xi(\tau)>m$ for all $x\in\R^N$, which contradicts the definition of~$m$.\par
Therefore,~(\ref{mzero}) holds. Now, as in the proof of Theorem~\ref{th1}, because of property~(\ref{hypf2}) at $z=m$, there exist~$R>0$ and~$\epsilon>0$ such that
$$\Delta(m+\epsilon\,\varphi_R)+f(m+\epsilon\,\varphi_R)\ge 0\hbox{ in }B(0,R)$$
and $m+\epsilon\,\varphi_R<u$ in $\overline{B(0,R)}$, where $\varphi_R$ solving~(\ref{phiR}) is the principal eigenfunction of the Dirichlet-Laplace operator in $B(0,R)$. Since $m+\epsilon\,\varphi_R=m$ on $\partial B(0,R)$ and $u>m$ in $\R^N$, the same sliding method as in the proof of Theorem~\ref{th1} implies that
$$m+\epsilon\,\varphi_R(x-y)\le u(x)\hbox{ for all }x\in\overline{B(y,R)}$$
and for all $y\in\R^N$. Therefore, $u\ge m+\epsilon\,\varphi_R(0)>m$ in $\R^N$, which contradicts the definition of~$m$.\par
As a conclusion, the assumption $u>m$ is impossible and, as already emphasized, the proof of Theorem~\ref{thspace} is thereby complete.\hfill$\Box$\break

The proof of Theorem~\ref{thhalfspace} also uses the sliding method and Theorem~\ref{thspace}, combined with limiting arguments as~$x_N\to+\infty$ and comparison with non-small subsolutions.\hfill\break

\noindent{\bf{Proof of Theorem~\ref{thhalfspace}.}} Let $u$ be a bounded nonnegative solution of~(\ref{eqhalf}) under assumption~(\ref{hypf2}). Since $f(0)\ge 0$, it follows from the strong maximum principle that either~$u=0$ in $\R^{N-1}\times\R_+$, or $u>0$ in $\R^{N-1}\times(0,+\infty)$. Let us then consider the second case. Since~$f(0)\ge 0$, it follows from Corollary~1.3 of Berestycki, Caffarelli and Nirenberg~\cite{bcn3} that~$u$ is increasing in~$x_N$. Denote
$$M=\sup_{\R^{N-1}\times[0,+\infty)}\,u.$$
There exists a sequence $(x_n)_{n\in\N}=(x_{1,n},\ldots,x_{N,n})_{n\in\N}$  in $\R^{N-1}\times\R_+$ such that $x_{N,n}\to+\infty$ and $u(x_n)\to M$ as $n\to+\infty$. From standard elliptic estimates, the functions
$$u_n(x)=u(x+x_n),$$
which satisfy the same equation as~$u$, converge in $C^{2,\beta}_{loc}(\R^N)$ for all $\beta\in[0,1)$, up to extraction of a subsequence, to a solution $u_{\infty}$ of~(\ref{eqRN}) such that $u_{\infty}(0)=M$. Theorem~\ref{thspace} implies that
$$u_{\infty}=M\hbox{ in }\R^N,$$
whence $f(M)=0$. It follows then from Berestycki, Caffarelli and Nirenberg~\cite{bcn1} (see also Theorem~1.4 in~\cite{bcn3}) that $u$ depends on $x_N$ only. In other words, the function $u(x)$ is equal to $V_M(x_N)$, which completes the proof of Theorem~\ref{thhalfspace}.\hfill$\Box$\break

Let us complete this section with the\hfill\break

\noindent{\bf{Proof of Theorem~\ref{thspacebis}.}} Assume that the function~$f$ satisfies~(\ref{hypf1}) and~(\ref{hypf3}) and let~$u$ be a bounded nonnegative solution of~(\ref{eqRN}). As already underlined, it first follows from the strong maximum principle that either $u\equiv 0$ in $\R^N$, or $u>0$ in $\R^N$. Let us then consider the second case. As in Step~1 of the proof of Theorem~\ref{th1} and applying the sliding method, one gets that
$$m=\inf_{\R^N}\,u>0.$$
Let $(x_n)_{n\in\N}$ be a sequence in $\R^N$ such that $u(x_n)\to m$ as $n\to+\infty$. Up to extraction of a subsequence, the functions
$$u_n(x)=u(x+x_n)$$
converge in $C^{2,\beta}_{loc}(\R^N)$ for all $\beta\in[0,1)$ to a classical bounded solution~$u_{\infty}$ of~(\ref{eqRN}) in~$\R^N$ such that $u_{\infty}\ge m$ in $\R^N$ and $u_{\infty}(0)=m$. Hence,
$$f(m)\le 0,$$
whence $m\ge\mu$ from~(\ref{hypf1}) and since $m>0$.\par
Call now
$$M=\sup_{\R^N}\,u.$$
If $M=m$, then $u$ is constant, which is the desired result. Assume now that $M>m$. The function
$$v=M-u$$
is a nonnegative bounded solution of
$$\Delta v+g(v)=0\ \hbox{ in }\R^N,$$
where the function~$g:[0,+\infty)\to\R$ is defined by
$$g(s)=\left\{\baa{ll}
-f(M-s) & \hbox{if }0\le s\le M-m,\vspace{3pt}\\
-f(m) & \hbox{if }s>M-m.\eaa\right.$$
Because of (\ref{hypf1}) and (\ref{hypf3}), the function~$g$ is Lipschitz-continuous and fulfills pro\-perty~(\ref{hypf2}). Theorem~\ref{thspace} applied to~$g$ and $v$ implies that the function~$v$ is actually constant. Hence,~$u$ is also constant, which actually shows that the assumption $M>m$ is impossible. As a conclusion, $M=m$ and~$u$ is then constant.\hfill$\Box$


\section{The dynamical systems' approach}\label{se4}

The goal of this section is to apply the dynamical systems' (shortly DS) approach to study the symmetrization and stabilization (as $x_1\rightarrow+\infty$) properties of the nonnegative solutions~(\ref{eq}) in
$$\Omega=(0,+\infty)\times\R^{N-2}\times(0,+\infty)$$
and (\ref{eqbis}) in
$$\Omega'=(0,+\infty)\times\R^{N-1}.$$
To this end we apply as aforementioned the DS approach. One of the main difficulties which arises in the dynamical study of~(\ref{eq}) in $\Omega$ or (\ref{eqbis}) in $\Omega'$ is the fact that the corresponding Cauchy problem is not well posed for~(\ref{eq}) in $\Omega$ and for (\ref{eqbis}) in $\Omega'$, and consequently the straightforward interpretation of~(\ref{eq}) and (\ref{eqbis}) as an evolution equation leads to semigroups of multivalued maps even in the case of cylindrical domains, see \cite{Babin}. The usage of multivalued maps can be overcome using the so-called trajectory dynamical approach (see \cite{Ball,cv,Sell} and the references therein). Under this approach, one fixes a signed direction~$\vec{l}$ in $\R^N$, which will play role of time. Then the space $K^+$ of all bounded nonnegative classical solutions of~(\ref{eq}) in $\Omega$ or (\ref{eqbis}) in $\Omega'$ (in the sense described in Section~\ref{intro}) is considered as a trajectory phase space for the semi-flow $(T^{\vec{l}}_h)_{h\in\R_+}$ of translations along the direction~$\vec{l}$ defined via
\begin{equation}
\big(T^{\vec{l}}_hu\big)(x)=u(x+h\vec{l}), \ \ h\in\R_+,\ u\in K^+.\label{4.1}
\end{equation}
In order the trajectory dynamical system $\big(T^{\vec{l}}_h,K^+\big)$ to be well defined, one needs the domains~$\Omega$ and~$\Omega'$ to be invariant with respect to positive translations along the~$\vec{l}$ directions, that is
\begin{align}
T^{\vec{l}}_h(\Omega)\subset\Omega\ \big(\hbox{resp. }T^{\vec{l}}_h(\Omega')\subset\Omega'\big),\ \ T^{\vec{l}}_h(x):=x+h\vec{l},\label{4.2}
\end{align}
for all $h\ge 0$. In our case, the $x_1$-axis will play the role of time, that is $\vec{l}=(1,0,\ldots 0,0)$. For the sake of simplicity of the notation, we then set
$$T^{\vec{l}}_h=T_h.$$\par
To apply the DS approach for our purposes, we apply the following Lemma~\ref{lem4.1} (see below), which also has an independent interest. For that purpose, let us introduce a few more notations. For any locally Lipschitz-continuous function~$f$ from $\R_+$ to $\R$, such that $f(0)\ge 0$, let $Z_f$ be defined by
\begin{equation}\label{defZf}
Z_f=\big\{z_0\in\R_{+}\ \big|\ f(z_0)=0\text{ and }F(z)<F(z_0)\text{ for }z\in[0,z_0)\big\},
\end{equation}
where
$$F(z)=\int^{z}_{0}f(\sigma)d\sigma.$$
The set $Z_f$ is then a subset of the set $E$ of zeroes of $f$, defined in~(\ref{defE}). Lastly, by~$R_f$ we denote the set of all bounded, nonnegative solutions $V\in C^2(\R_+)$ of
\begin{equation}\left\{\baa{l}
        V''(\xi)+f(V(\xi))=0\ \text{for all } \xi\geq 0,\vspace{3pt}\\
        V(0)=0,\ \ V\geq 0,\ \ V\ \text{is bounded}.\eaa\right.\label{4.3}
\end{equation}

\begin{lem}\label{lem4.1} Let $f$ be a locally Lipschitz-continuous function from $\R_+$ to $\R$, such that $f(0)\ge 0$. Then the set $R_f$ is homeomorphic to $Z_f$ and as a consequence is totally disconnected.
\end{lem}

\noindent\textbf{Proof.} Since $f(0)\ge 0$, it follows from elementary arguments that every nonnegative bounded solution of (\ref{4.3}) is monotonically nondecreasing, that is $V(\xi_1)\geq V(\xi_2)$ if $\xi_1\geq\xi_2$. Consequently, the following limit
\begin{align}
z_0=z_0(V)=\underset{\xi\rightarrow+\infty}{\lim}V(\xi)\label{4.4}
\end{align}
exists and necessarily $f(z_0)=0$ and $0\leq V(z)\leq z_0$. Moreover, it is well known that either $V(\xi)=0$ for all $\xi\ge 0$, or $V'(\xi)>0$ for all $\xi\geq 0$.\par
Multiplying equation~(\ref{4.3}) by $V'$ and integrating over $[0,\xi]$ we obtain the explicit expression for the derivative $V'(\xi)$:
\begin{align}
V'(\xi)^2=-2F(V(\xi))+V'(0)^2\ \hbox{ for all }\xi\ge 0.\label{4.5}
\end{align}
Passing to the limit $\xi\rightarrow+\infty$ in~(\ref{4.5}) and taking into account~(\ref{4.4}), it follows that $V'(0)^2=2F(z_0)$, whence $F(z_0)\ge 0$. Therefore we obtain the following equation for $V(\xi)$:
\begin{align}
V'(\xi)^2=2\big(F(z_0)-F(V(\xi))\big)\ \hbox{ for all }\xi\ge 0.\label{4.6}
\end{align}\par
Assume now that $F(z_0)>0$, whence $z_0>0$ and then $V'>0$ in $\R_+$. Then the solution~$V_{z_0}$ of~(\ref{4.3}) that satisfies~(\ref{4.4}) and~(\ref{4.6}) exists if and only if $F(z_0)-F(z)>0$ for every $z\in[0,z_0)$. Moreover, such a solution is unique because $V_{z_0}$ satisfies~(\ref{4.3}) with the initial conditions
\begin{align}
V_{z_0}(0)=0, V'_{z_0}(0)=\sqrt{2F(z_0)}.\label{4.7}
\end{align}
If $z_0=0$, then $V'(0)=0$ and $V(\xi)=0$ for all $\xi\in\R_+$, whence $f(0)=0$. In this case, we define $V_0=0$ in $\R_+$.\par
Next we show that $Z_f$ defined by (\ref{defZf}) is totally disconnected in $\R$. Indeed, otherwise it should contain a segment $[a,b]$ with $0\le a<b$. Then, $f(z_0)\equiv 0$ for $z_0\in[a,b]$ and consequently $F(z_0)=F(b)$ for every $z_0\in[a,b]$, which evidently leads to a contradiction.\par
To prove the disconnectedness of the set $R_f$ it is then sufficent to show that there exists a homeomorphism
$$\tau: \big(Z_f,\R\big)\rightarrow\big(R_f,C^2_{loc}(\R_+)\big).$$
To do so, observe that (\ref{4.7}) defines a homeomorphism between $Z_f$ and the set
$$R_f(0)=\big\{(0,V'(0))\ \big|\ V\in R_f\big\}$$
of values at $\xi=0$ of the functions from $R_f$. Recall that $R_f$ consists of the solutions of the second order ODE~(\ref{4.3}) and, consequently, thanks to the classical Cauchy-Lipschitz theorem on continuous dependence of solutions of ODE's, the set $R_f$ is homeomorphic to~$R_f(0)$ and this homeomorphism is given by $V\mapsto(0,V'(0))$. As a conclusion, the set~$R_f$ is totally disconnected.\hfill$\Box$

\begin{rem}\label{rem4.1} {\rm Note that, although for generic functions $f$'s the set $Z_f \cong R_f$ is finite, this set may be even uncountable for some very special choices of nonlinearities $f$. One of the simplest examples of such a function $f$ is the following one:
\begin{align}
f(z):=dist (z,\cal C)\label{4.9}
\end{align}
for all $z\in\R_+$, where $\cal C$ is a standard Cantor set on $[0,1]$ and $dist(z,\cal C)$ means a distance from $z$ to $\cal C$. Indeed, it is easy to verify that for this case $Z_f=\cal C$ and consequently $R_f$ consists of a continuum of elements.}
\end{rem}

Below we state the main results of this Section~\ref{se4}, that is Theorems \ref{the4.1}-\ref{the4.4}, which are obtained by the dynamical systems' approach. To this end we define a class of functions $K^+$ to which the solutions of~(\ref{eq}) as well as (\ref{eqbis}) belong to. Namely, a bounded nonnegative solution of~(\ref{eq}) (resp. (\ref{eqbis})) is understood to be a solution $u$ of class $C^2(\Omega)$ (resp. $C^2(\Omega')$) and continuous on $\overline{\Omega}\,\backslash\,\{0\}\!\times\!\R^{N-2}\!\times\!\{0\}$ (resp. on $\overline{\Omega'}$). The set $K^+$ is endowed with the local topology according to the embedding of $K^+$ in $C^{2,\beta}_{loc}(\overline{\Omega}\backslash\{x_1=0\})$ (resp. $C^{2,\beta}_{loc}(\overline{\Omega'}\backslash\{x_1=0\})$) for all $\beta\in[0,1)$. Actually, as already emphasized in Section~\ref{intro}, all solutions $u\in K^+$ are automatically in $C^{2,\beta}_b([\epsilon,+\infty)\times\R^{N-2}\times[0,+\infty))$ (resp. $C^{2,\beta}_b([\epsilon,+\infty)\times\R^{N-1})$) for all $\beta\in[0,1)$ and for all $\epsilon>0$, where we refer to~(\ref{defC2betab}) for the definition of the sets~$C^{2,\beta}_b(F)$.\par
The first two theorems are concerned with the case of functions $f$ fulfilling the condition~(\ref{hypf1}).

\begin{theo} \label{the4.1} Let $N$ be any integer such that $N\ge 2$ and let $f$ be a locally Lipschitz-continuous function from $\R_+$ to $\R$, satisfying~$(\ref{hypf1})$. Then the trajectory dynamical system~$\left(T_h,K^+\right)$ associated to~$(\ref{eq})$ possesses a global attractor $A_{tr}$ in $K^+$ which is bounded in~$C^{2,\beta}_{b}(\overline{\Omega})$ and then compact in $C^{2,\beta}_{loc}(\overline{\Omega})$ for all $\beta\in[0,1)$. Moreover~$A_{tr}$ has the following structure
$$A_{tr}=\Pi_{\overline{\Omega}}K^+\big(\widetilde{\Omega}\big)$$
where $\widetilde{\Omega}=\R^N_+$, $K^+\big(\widetilde{\Omega}\big)$ is the set of all bounded nonnegative solutions of~$(\ref{eqhalf})$ in $\widetilde{\Omega}=\R^N_+$, and $\Pi_{\overline{\Omega}}$ denotes the restriction to $\overline{\Omega}$. Hence,
$$A_{tr}\subset\big\{x\mapsto 0,\ x\mapsto V(x_N)\big\}$$
and $A_{tr}=\{x\mapsto V(x_N)\}$ if $f(0)>0$, where $V$ is the unique solution of~$(\ref{eq1d})$. Lastly, for any bounded nonnegative solution~$u$ of~$(\ref{eq})$ in~$\Omega$, the functions $T_hu$ converge as $h\rightarrow+\infty$ in $C^{2,\beta}_{loc}(\overline{\Omega})$ for all $\beta\in[0,1)$ either to $0$ or to $x\mapsto V(x_N)$, and they do converge to the function $x\mapsto V(x_N)$ if $f(0)>0$.
\end{theo}

The next theorem deals with the analysis of the asymptotic behavior as $x_1\rightarrow+\infty$ of the nonnegative bounded classical solutions $u\in K^+$ of equation~(\ref{eqbis}) in the half-space $\Omega'=(0,+\infty)\times\R^{N-1}$. For any $M\ge 0$, we define
$$K^+_M=K^+\,\cap\,\big\{0\leq u\leq M\big\}.$$

\begin{theo}\label{the4.2} Let $N$ be any integer such that $N\geq 2$ and assume that, in addition to~$(\ref{hypf1})$, the given locally Lipschitz-continuous function $f$ from $\R_+$ to $\R$ satisfies~$(\ref{hypf3})$. Then, for every $M\ge\mu$, the trajectory dynamical system $\left(T_h,K^+_M\right)$ associated to~$(\ref{eqbis})$ possesses a global attractor $A_{tr}$, which is bounded in $C^{2,\beta}_{b}(\overline{\Omega'})$ and then compact in $C^{2,\beta}_{loc}(\overline{\Omega'})$ for all $\beta\in[0,1)$, and satisfies
\begin{align}
A_{tr}=\Pi_{\overline{\Omega'}}K^+_M\big(\widetilde{\Omega'}\big)\label{projection}
\end{align}
where $\widetilde{\Omega'}=\R^N$, $K^+_M\big(\widetilde{\Omega'}\big)=K^+\big(\widetilde{\Omega'}\big)\cap\big\{0\le u\le M\big\}$, $K^+\big(\widetilde{\Omega'}\big)$ is the set of all bounded nonnegative solutions of~$(\ref{eqRN})$ in $\widetilde{\Omega'}=\R^N$, and $\Pi_{\overline{\Omega'}}$ denotes the restriction to $\overline{\Omega'}$. Hence,
$$A_{tr}=\big\{z\in[0,M]\ \big|\ f(z)=0\big\}=E\cap[0,M].$$
Lastly, for any bounded nonnegative solution $u$ of~$(\ref{eqbis})$ in $\Omega'$ such that $0\le u\le M$, the functions~$T_hu$ converge as $h\to+\infty$ in $C^{2,\beta}_{loc}(\overline{\Omega'})$ for all $\beta\in[0,1)$ to some $z\in E\cap[0,M]$ which is uniquely defined by~$u$. 
\end{theo}

The next two theorems, which are concerned with the case of functions $f$ fulfilling the condition~(\ref{hypf2}), are based on the new Liouville type Theorems~\ref{thspace} and \ref{thhalfspace} which were already stated in Section~\ref{sec2}.

\begin{theo} \label{the4.3} Let $N$ be any integer such that $N\geq 2$, let $f$ be any locally Lipschitz-continuous function from $\R_+$ to $\R$ satisfying~$(\ref{hypf2})$ and assume that, for problem~$(\ref{eq})$ in the quarter-space~$\Omega$, the set $K^+$ is not empty. Then, for every sufficiently large $M\ge 0$, the trajectory dynamical system $\big(T_h,K^+_M\big)$ associated to~$(\ref{eq})$ possesses a global attractor~$A_{tr}$, which is bounded in~$C^{2,\beta}_{b}(\overline{\Omega})$ and then compact in $C^{2,\beta}_{loc}(\overline{\Omega})$ for all $\beta\in[0,1)$ and has the following structure
\be\label{projectionbis}
A_{tr}=\Pi_{\overline{\Omega}}K^+_M\big(\widetilde{\Omega}\big)
\ee
where $K^+_M\big(\widetilde{\Omega}\big)=K^+\big(\widetilde{\Omega}\big)\cap\big\{0\leq u\leq M\big\}$. Hence,
$$A_{tr}=\big\{x\mapsto V_z(x_N)\ |\ z\in[0,M],\ f(z)=0\big\}.$$
Lastly, for any bounded nonnegative solution~$u$ of~$(\ref{eq})$ in $\Omega$ such that $0\le u\le M$, the functions~$T_hu$ converge as $h\to+\infty$ in $C^{2,\beta}_{loc}(\overline{\Omega})$ for all $\beta\in[0,1)$ to some function $x\mapsto V_z(x_N)$, where $z\in E\cap[0,M]$ is uniquely defined by~$u$ and $E$ denotes the set of zeroes of the function~$f$.
\end{theo}

Analogously to Theorem~\ref{the4.2} we have the following Theorem~\ref{the4.4} in the case of the half-space $\Omega'=(0,+\infty)\times\R^{N-1}$.

\begin{theo} \label{the4.4} Let $N$ be any integer such that $N\geq 2$, let $f$ be any locally Lipschitz-continuous function from $\R_+$ to $\R$ satisfying~$(\ref{hypf2})$, and assume that, for problem~$(\ref{eqbis})$ in the half-space~$\Omega'$, the set $K^+$ is not empty. Then, for every sufficiently large $M\ge 0$, the trajectory dynamical system $\big(T_h,K^+_M\big)$ associated to~$(\ref{eqbis})$ possesses a global attractor~$A_{tr}$, which is bounded in $C^{2,\beta}_{b}(\overline{\Omega'})$ and then compact in $C^{2,\beta}_{loc}(\overline{\Omega'})$ for all $\beta\in[0,1)$, and satisfies~$(\ref{projection})$. Lastly, for any bounded nonnegative solution~$u$ of~$(\ref{eqbis})$ in $\Omega'$ such that $0\le u\le M$, the functions~$T_hu$ converge as $h\to+\infty$ in $C^{2,\beta}_{loc}(\overline{\Omega'})$ for all $\beta\in[0,1)$ to some $z\in E\cap[0,M]$ which is uniquely defined by~$u$.
\end{theo}

In what follows we prove Theorem~\ref{the4.1}. A proof of Theorems~\ref{the4.2}-\ref{the4.4} can be done in the same manner as in Theorem~\ref{the4.1} with some minor modifications (see also Remark~\ref{rem4.3}). In particular, in Theorems~\ref{the4.3} and~\ref{the4.4}, one can take~$M$ as any nonnegative real number such that $M\ge\|U\|_{\infty}$, where~$U$ is any element in~$K^+$. We leave the details to the reader.\hfill\break

\noindent\textbf{Proof of Theorem~\ref{the4.1}.} Let $K^+$ be the set of all bounded nonnegative solutions of~(\ref{eq}) in~$\Omega$. Due to the assumptions (\ref{hypf1}), the set $K^+$ is not empty (the function $x\mapsto V(x_N)$, where~$V$ is the unique solution of~(\ref{eq1d}), belongs to $K^+$) and due to the translation invariance of~(\ref{eq}), it follows that $T_h:K^+\rightarrow K^+$ is well defined for all $h\ge 0$, where $(T_h u)(x_1,x',x_N):=u(x_1+h,x',x_N)$.\par
To show that $\left(T_h,K^+\right)$ possesses a global attractor, it suffices to show (see~\cite{Babin} and the references therein) that
\begin{itemize}
\item for any fixed $h>0$, $T_h$ is a continuous map in $K^+$ (we recall that $K^+$ is endowed with local topology according to the embedding of~$K^+$ in $C^{2,\beta}_{loc}(\overline{\Omega}\backslash\{x_1=0\})$) for all $\beta\in[0,1)$;
\item the semi-flow $(T_h)_{h\ge 0}$ possesses a compact attracting (absorbing) set in $C^{2,\beta}_{loc}(\overline{\Omega})$, which is even bounded in $C^{2,\beta}_b(\overline{\Omega})$, for all $\beta\in[0,1)$.
\end{itemize}
Note that, the continuity of $T_h$ in $C^{2,\beta}_{loc}(\overline{\Omega}\backslash\{x_1=0\})$ is obvious, because the shift operator is continuous in this topology, as well as its restriction to~$K^+$. As for the existence of compact attracting (absorbing) set for the semi-flow $(T_h)_{h\ge 0}$, it follows from the fact that the set of all bounded nonnegative solutions of~(\ref{eqhalf}) in $\cup_{h\ge 0}T_{-h}(\overline{\Omega})=\widetilde{\Omega}=\R^N_+$ under the assumption~(\ref{hypf1}) on $f$ is uniformly bounded. Indeed, as already recalled in Section~\ref{sec21} and according to a result of~\cite{bcn2}, under the assumption (\ref{hypf1}), any bounded solution of~(\ref{eqhalf}) in~$\R^N_+$ which is positive in $\R^{N-1}\times(0,+\infty)$ has one-dimensional symmetry, that is
$$u(x_1,x',x_N)=V(x_N)$$
where $0\leq V<\mu$ is the unique solution of~(\ref{eq1d}). On the other hand, since $f(0)\ge 0$, any bounded nonnegative solution of~(\ref{eqhalf}) is either positive in $\R^{N-1}\times(0,+\infty)$, or identically~$0$ in~$\R^N_+$, and it cannot be $0$ if $f(0)>0$. Thus, the set $K^+\big(\widetilde{\Omega}\big)$ of all bounded nonnegative solutions of (\ref{eqhalf}) in $\widetilde{\Omega}=\R^N_+$ is bounded in $L^{\infty}(\widetilde{\Omega})$, namely
$$\sup_{u\in K^+(\widetilde{\Omega})}\|u\|_{\infty}\le\mu.$$
Then the existence of a compact absorbing set for $(T_h)_{h\ge 0}$ in $C^{2,\beta}_{loc}(\overline{\Omega})$, which is even bounded in $C^{2,\beta}_b(\overline{\Omega})$, for all $\beta\in[0,1)$ is a consequence of the uniform boundedness of all solutions of (\ref{eqhalf}) in $\widetilde{\Omega}=\R^N_+$ and of standard elliptic estimates. Hence the semigroup $\big(T_h,K^+\big)$ possesses a global attractor $A_{tr}$ in $K^+$ which is bounded in $C^{2,\beta}_b(\overline{\Omega})$ and compact in $C^{2,\beta}_{loc}(\overline{\Omega})$ for all $\beta\in[0,1)$.\par
To prove the convergence part of Theorem~\ref{the4.1}, as we will see below, it is sufficient to show that $A_{tr}=\Pi_{\overline{\Omega}}K^+(\widetilde{\Omega})$. Assuming for a moment that this representation is true, we obtain from the previous considerations that
\be\label{subset}
A_{tr}\subset\big\{x\mapsto 0,\ x\mapsto V(x_N)\big\},
\ee
and $A_{tr}$ is then equal to the singleton $\big\{x\mapsto V(x_N)\big\}$ if $f(0)>0$. Hence, for any bounded nonnegative solution $u$ of (\ref{eq}), since $\big\{T_hu,h\geq 1\big\}$ is bounded in $C^{2,\beta}_b(\overline{\Omega})$ and compact in $C^{2,\beta}_{loc}(\overline{\Omega})$ for all $\beta\in[0,1)$, the $\omega$-limit set $\omega(u)$ of $u$ is not empty and it is an invariant and connected subset of $A_{tr}$. Since $A_{tr}$ is totally disconnected,\footnote{This property is obvious here due to (\ref{subset}). See Remark~\ref{rem4.3} for a comment about the other situations, corresponding to Theorems~\ref{the4.2},~\ref{the4.3} and~\ref{the4.4}.} it follows that either $\omega(u)=\big\{x\mapsto 0\big\}$ or $\omega(u)=\big\{x\mapsto V(x_N)\big\}$, the latter being necessarily true if $f(0)>0$.\par
To complete the proof of Theorem~\ref{the4.1} it remains to show that $A_{tr}=\Pi_{\overline{\Omega}}K^+(\widetilde{\Omega})$. First we prove that $\Pi_{\overline{\Omega}}\hat{u}\in A_{tr}$ for any bounded nonnegative solution $\hat{u}$ of~(\ref{eqhalf}) in $\widetilde{\Omega}=\R^N_+$, that is $\hat{u}\in K^+(\widetilde{\Omega})$. Indeed, for such a $\hat{u}$, the family $(\Pi_{\overline{\Omega}}(T_{-h}\hat{u}))_{h\ge 0}$ is uniformly bounded in~$C^{2,\beta}_b(\overline{\Omega})$ and compact in $C^{2,\beta}_{loc}(\overline{\Omega})$ for all $\beta\in[0,1)$ and, according to definition of the attractor, there holds
$$T_h\Pi_{\overline{\Omega}}(T_{-h}\hat{u})\longrightarrow A_{tr}\text{ in }C^{2,\beta}_{loc}(\overline{\Omega})\text{ as }h \rightarrow+\infty,$$
for all $\beta\in[0,1)$. On the other hand, $T_h\Pi_{\overline{\Omega}}(T_{-h}\hat{u})=\Pi_{\overline{\Omega}}\hat{u}$. Hence $\Pi_{\overline{\Omega}}\hat{u}\in A_{tr}$.\par
Next we prove the reverse inclusion. To this end, let us recall that $\left(T_h,K^+\right)$ possesses an absorbing set which is bounded in $C^{2,\beta}_b(\overline{\Omega})$ and then compact in $C^{2,\beta}_{loc}(\overline{\Omega})$ for all $\beta\in[0,1)$, say $\B_*\subset K^+$, and, as a consequence,
$$A_{tr}=\omega(\B_*)=\underset{h\geq 0}{\bigcap}\left[\underset{s\geq h}{\bigcup}T_s \B_*\right],$$
where $[$ $]$ means the closure in $C^{2,\beta}_{loc}(\overline{\Omega})$ (see~\cite{Babin,cv} and the references therein). Let now $u\in A_{tr}$. The property $A_{tr}=\omega(\B_*)$ implies that there exist an increasing sequence $(h_k)_{k\in\N}\rightarrow+\infty$ and a sequence of solutions $(u_k)_{k\in\N}$ in $\B_*$, such that
\begin{align}
u=\lim_{k\to+\infty}T_{h_k}u_k\label{*}
\end{align}
in $C^{2,\beta}_{loc}(\overline{\Omega})$ for all $\beta\in[0,1)$. Note that the solution $T_{h_k}u_k$ is defined not only in $\Omega$, but also in the domain $(-h_k,+\infty)\times\R^{N-2}\times\R_+$, and that
\begin{align}
\sup_{k\in\N}\ \|T_{h_k}u_k\|_{C^{2,\beta}_{b}\left([-h_k+\epsilon,\infty)\times\R^{N-2}\times\R_+\right)}<+\infty\label{**}
\end{align}
for all $\epsilon>0$ and $\beta\in[0,1)$, from standard elliptic estimates. Consequently, for every~$k_0\in\N$ and $\beta\in[0,1)$, the sequence $(T_{h_k}u_k)_{k>k_0}$ is precompact in $C^{2,\beta}_{loc}\left([-h_{k_0},\infty)\times\R^{N-2}\times\R_+\right)$. Taking a subsequence, if necessary, and using Cantor's diagonal procedure and the fact $h_k\rightarrow\infty$, we can say that this sequence converges to $\hat{u}\in C^{2,\beta}_{loc}(\R^N_+)$ in the spaces $C^{2,\beta}_{loc}\left([-h_{k_0},\infty)\times\R^{N-2}\times\R_+\right)$ for every $k_0\in\N$ and for every $\beta\in[0,1)$. Then~(\ref{**}) implies that $\hat{u}\in C^{2,\beta}_{b}(\R^N_+)$ for every $\beta\in[0,1)$. Lastly, the functions~$T_{h_k}u_k$ are nonnegative solutions of~(\ref{eq}) in $(-h_k,+\infty)\times\R^{N-2}\times\R_+$ and by letting $k\rightarrow+\infty$, we easily obtain that~$\hat{u}$ is a bounded nonnegative solution of~(\ref{eqhalf}) in $\widetilde{\Omega}=\R^N_+$. Finally, formula~(\ref{*}) implies that
$$\Pi_{\overline{\Omega}}\hat{u}=u.$$
Thus $u\in\Pi_{\overline{\Omega}}K^+(\widetilde{\Omega})$ and the representation formula $A_{tr}=\Pi_{\overline{\Omega}}K^+(\widetilde{\Omega})$ is proved. The proof of Theorem~\ref{the4.1} is thereby complete.\hfill$\Box$

\begin{rem} \label{rem4.2} {\rm As far as Theorems~\ref{th1} and~\ref{th2} on the one hand, and Theorems~\ref{the4.1} and~\ref{the4.3} on the other hand, are concerned, we especially emphasize that the DS approach simplified in a very elegant way most of the computations regarding the asymptotic behavior of the solutions of~(\ref{eq}) as $x_1\rightarrow+\infty$. However, the DS approach and Theorem~\ref{the4.1} (resp. Theorem~\ref{the4.3}) do not provide as in PDE approach the fact that only the limiting profile~$V(x_N)$ (resp. $V_z(x_N)$ for some $z\in E\backslash\{0\}$) is selected, even if $f(0)=0$, as soon as $u_0\not\equiv0$ on~$\{0\}\times\R^{N-2}\times(0,+\infty)$. In the PDE proof of Theorems~\ref{th1} and~\ref{th2}, it is indeed shown that the condition $u_0\not\equiv0$ implies that $u$ is separated from $0$ for large enough $x_1$ and $x_N$. This property is not shown in the DS proof of Theorems~\ref{the4.1} and~\ref{the4.3}. Similar comments also hold for Theorems~\ref{th3},~\ref{th4},~\ref{the4.2} and~\ref{the4.4}, where the PDE proof provides the convergence to a {\it non-zero} zero of~$f$, what the DS proof does not. Lastly, in some of the results obtained through the PDE approach, the convergence of the solutions as $x_1\to+\infty$ is proved to be uniform with respect to the variables $(x',x_N)$, while the DS approach only provides local convergence, due to the necessity of using the local topology to get the existence of a global attractor.}
\end{rem}

\begin{rem} \label{rem4.3} {\rm In Theorem~\ref{the4.3} (resp. Theorem~\ref{the4.4}) under assumption~(\ref{hypf2}), the phase space~$K^+_M$, which is invariant under the semigroup $T_h$, is not empty for any sufficiently large~$M$, because $K^+$ is assumed to be not empty. Then, in the same manner as in the proof of Theorem~\ref{the4.1}, using both the representation formula $A_{tr}=\Pi_{\overline{\Omega}}K^+_M(\widetilde{\Omega})$ (resp. $A_{tr}=\Pi_{\overline{\Omega'}}K^+_M(\R^N)$), the Liouville theorems of Section~\ref{sec2} and Lemma~\ref{lem4.1} one obtains the desired conclusions. In particular, for Theorem~\ref{the4.3} (resp. Theorem~\ref{the4.4}) about problem~(\ref{eq}) in~$\Omega$ (resp.~(\ref{eqbis}) in $\Omega'$), the total disconnectedness of $A_{tr}$ follows from the representation formula~(\ref{projectionbis}) (resp.~(\ref{projection})), from Theorem~\ref{thhalfspace} (resp. Theorem~\ref{thspace}) and from Lemma~\ref{lem4.1} with, here, $Z_f=E$ (resp. condition~(\ref{hypf2}) again). Note that for Theorem~\ref{the4.2} in the case of assumptions~(\ref{hypf1}) and~(\ref{hypf3}), the total disconnectedness of~$A_{tr}$ follows from Theorem~\ref{thspacebis} and assumption~(\ref{hypf3}) again.}
\end{rem}

\begin{rem} \label{rem4.4} {\rm Note that neither the concrete choice of the domain  $\Omega$ (or $\Omega'$) nor the concrete choice of the ``time" direction $x_1$ are essential for the use of the trajectory dynamical system' approach. Indeed, let us replace the ``time" direction $x_1$  by any fixed direction $\vec l\in \R^N$ and correspondingly $T^{\vec{l}}_hu=u(\cdot+h\vec{l}\,)$ for $h\in\R_+$ and $u\in K^+$. Then the above construction seems to be applicable if the domain $\Omega$ satisfies the following assumptions:
\begin{itemize}
\item $T_h\Omega\subset\Omega$ (this is necessary in order to define the restriction of $T_h$ to the trajectory phase space $K^+$ or $K^+_M$).
\item $\cup_{h\ge 0}T_{-h}\overline{\Omega}=\R^N_+$, or $\R^N$ (this is required in order to obtain representation formulas of the type $A_{tr}=\Pi_{\overline{\Omega}}K^+(\R^N_+)$ or $A_{tr}=\Pi_{\overline{\Omega}}K^+(\R^N)$, with possibly $K^+_M$ instead of~$K^+$).
\end{itemize}}
\end{rem}


\end{document}